\documentclass{amsart}

\usepackage{graphicx}
\usepackage{hyperref}
\usepackage{multicol}
\usepackage{multirow}
\usepackage{xcolor}
\colorlet{RED}{red}

\newtheorem{theorem}{Theorem}[section]

\newtheorem{conjecture}[theorem]{Conjecture}

\theoremstyle{definition}

\theoremstyle{remark}

\def\Z{\mbox{$\mathbb Z$}}

\numberwithin{equation}{section}

\begin{document}

\title{The Exit Distribution for Smart Kinetic Walk with Symmetric and Asymmetric Transition Probability}

\author{Yan Dai}
\email{ydai@math.arizona.edu}
\address{Department of Mathematics, University of Arizona, Tucson, AZ 85719}

\begin{abstract}
It has been proved that the distribution of the point where the Smart Kinetic Walk (SKW) exits a domain converges in distribution to harmonic measure on the hexagonal lattice. For other lattices, it is believed that this result still holds, and there is good numerical evidence to support this conjecture. Here we examine the effect of the symmetry  and asymmetry of the transition probability on each step of the SKW on the square lattice and test if the exit distribution converges in distribution to harmonic measure as well. From our simulations, the limiting exit distribution of the SKW with a non-uniform but symmetric transition probability as the lattice spacing goes to zero is the harmonic measure. This result does not hold for asymmetric transition probability. We are also interested in the difference between the SKW with symmetric transition probability exit distribution and harmonic measure. Our simulations provide strong support for a explicit conjecture about this first order difference. The explicit formula for the conjecture will be given below.
\end{abstract}

\maketitle

\section{Introduction}
The ordinary random walk  has been well studied, and it is a widely known theorem that the exit distribution of the ordinary random walk converges to harmonic measure. A proof of this theorem can be found in \cite{LL}. There are some other interesting random walks, such as random walks that have some self-avoiding constraint. One such model is the self-avoiding walk (SAW) such that the walk visits each site at most once, and all walks with the same length have the same probability \cite{MS}. The scaling limit of the SAW in a simply connected domain is conjectured to be Schramm-Loewner evolution (SLE) with parameter $\kappa$=8/3 \cite{LSW}. This result has a good numerical support \cite{K02, K04}.  In this paper, we study another type of random walk called smart kinetic walk (SKW), which never intersects itself and never steps into a site that would lead to its being trapped. It first appeared in the physics literature in mid 1980s under some different names. In \cite{KL}, Kremer and Lyklema called it the infinitely growing SAW. Their model is the model that we consider in the full plane. Weinrib and Trugman have independently introduced the model and called it the smart kinetic walk \cite{WT}.  Ziff, Cummings, and Stell introduced a two-dimensional random walk that can be used to generate the perimeter of percolation cluster for critical site percolation on the square lattice \cite{ZCS}. A walk generated by their random walk algorithm is different from a walk under the SKW model since the occupied site can be visited more than once. Grassberger studied a random walk model on the square lattice which only turns left or right with equal probability, straight continuations are forbidden. This model generates the perimeters for bond percolation on the square lattice \cite{G}. It is not quite the same as the SKW model in this paper. In our model, no site is visited more than once. We can refer to it as the site avoiding walk. Grassberger's model is the bond avoiding walk, which forbids the walk to visit a bond more than once. Walks with loops are allowed under his model. Note the SKW model can be defined in any lattice and any dimension. We are interested in the two-dimensional walk on the square lattice. 

Let $D$ be a bounded domain with smooth boundary containing the origin. To generate a SKW $\omega=(\omega(0), \omega(1), \cdots,\omega(n), \cdots) $ with each $\omega(i)\in \Z^2$ in $D$, we start the walk at the origin, i.e., $\omega(0)=0$. In the first step, the walk can choose uniformly from the nearest neighbors of the origin with equal probability. The walk stops as soon as it reaches the boundary of the domain $D$. If a site has been visited before or is outside of the domain, we call it an occupied site. A trapping site is a site that would lead it to get trapped, i.e., there is no path from the site to the boundary of the domain through unoccupied sites. If a site is neither an occupied site nor a trapping site, we call it as an allowable site. 

\begin{figure}[h!]
  \centering
      \includegraphics[width=0.65\textwidth]{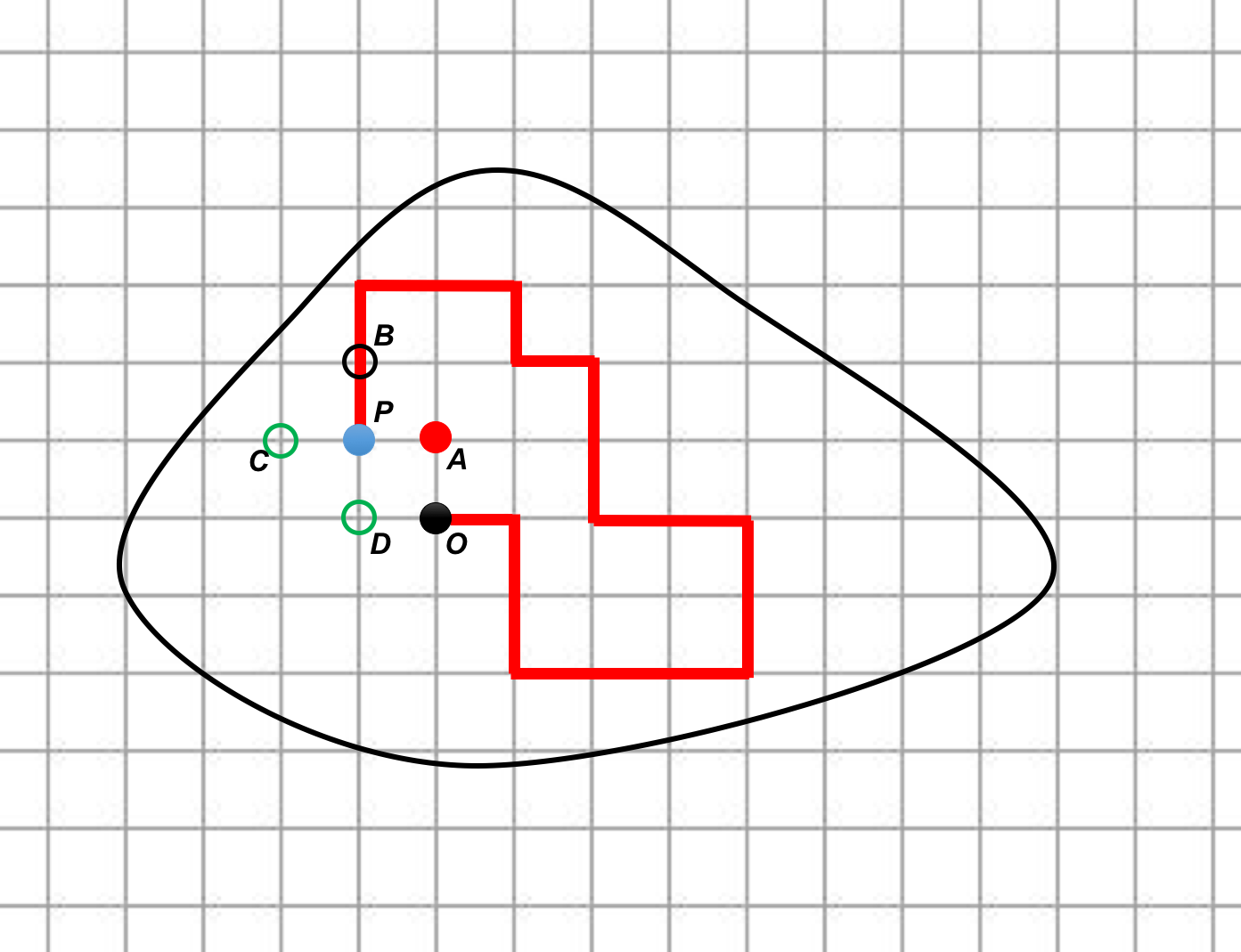}
  \caption{An example of SKW illustrating occupied site, allowable site and trapping site.}
  \label{skw_def}
\end{figure}

Consider the allowable nearest neighbors of $\omega(n)$ and assign them with equal probability. We pick one randomly and assign it to be $\omega(n+1)$. Figure \ref{skw_def}  illustrates a SKW starting at the origin $O$ and walking up to the site $P$. The nearest neighbor $B$ is occupied. Site $A$ is a trapping site since it has no allowable path to the boundary of the domain. The walker can choose $C$ or $D$ for the next step each with probability 1/2. 

On the hexagonal lattice, the SKW in a bounded domain is equivalent to a percolation interface \cite{CN, W}. It has been proved that the exit distribution of the SKW converges in distribution to harmonic measure on the hexagonal lattice \cite{Jiang}.  In \cite{tgk} the usual SKW model has been modified by averaging over the orientation of the lattice with respect to the domain, and Monte Carlo simulations have been carried out for this modified model on three lattices: triangular, square and hexagonal. From these simulations, this result still holds. Moreover, the simulations lead to a conjecture that the difference between harmonic measure and the exit distribution of the rotationally averaged SKW model, to first order in the lattice spacing $\delta$, is given by $\delta c \rho_D(z)|dz|$, where $c$ is a constant that only depends on the lattice and the model, and the density function $\rho_D$ depends only on the domain. See \cite{tgk} for details. 

In this paper, we examine the effect of the symmetry and asymmetry  of the transition probability on each step of the SKW on the square lattice and test if the exit distribution converges weakly to harmonic measure. We are also interested in the leading order term in the difference between the SKW with symmetric transition probability exit distribution and harmonic measure, which we refer to as the first order correction. In our model, we also average over rotations of the orientation of the lattice with respect to the domain $D$. Comparing with some other models with fixed lattice orientation, since the scaling limit is rotationally invariant, the exit distribution does not depend on the orientation of the lattice. However, there is no reason to expect that the first order correction is independent to the lattice orientation. Therefore, we are not studying the first order correction with the fixed orientation of the lattice. Instead, we are studying rotationally averaged SKW model here.

\begin{figure}[h!]
  \centering
      \includegraphics[width=0.8\textwidth]{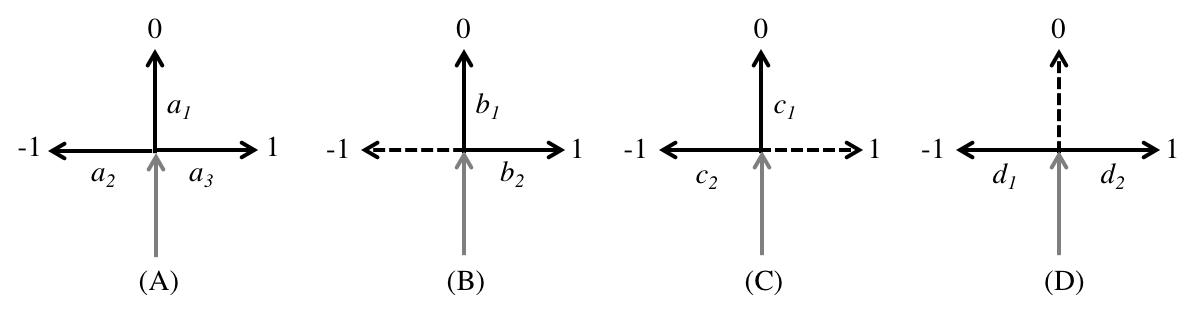}
  \caption{The transition probability in the SKW model.}
  \label{transitionprob}
\end{figure}

Note that in the original SKW model, the transition probability at each step of the walk is defined as $p$ = 1/(number of the allowable sites).  The walker starts at the origin and chooses one of the sites from four of the nearest neighbors each with probability 1/4. After that, at each step, the walker searches for allowable sites. Figure \ref{transitionprob} shows possible allowable sites in different situations. The dashed line in the figure indicates that the site is blocked and the grey line indicates the direction of the previous step. In the original SKW model, if all three neighbors are allowable, then the transition probability is 1/3. The variables $a_1, a_2$ and $a_3$ are all equal to 1/3 in figure \ref{transitionprob} (A) . If one of the sites (left, right or front) is blocked, which corresponds to (B), (C) and (D) in the figure \ref{transitionprob} respectively, then the transition probability is 1/2.  This says the variables in (B), (C) and (D) are all equal to 1/2.  We give a name to each case in figure \ref{transitionprob}:  (A) as the nblock case, (B) as the left-blocked case, (C) as the right-blocked case and (D) as the front-blocked case. We have not bothered to draw the trivial case for which there is only one allowable site. In the original model, the transition probability is symmetric and uniform. In this paper, we consider the transition probabilities that are not uniform any more. 

We say that the transition probability in the SKW is symmetric if it satisfies the following three conditions:
\begin{itemize}
\item[(i)] $a_2=a_3$,
\item[(ii)] $b_1=c_1$, which implies $b_2=c_2$,
\item[(iii)] $d_1=d_2=1/2$.
\end{itemize}

Given a domain $D$, there are different ways to define the point where the walk first exits the domain since it may not be exactly on the boundary of the domain $D$. In this paper, we use orthogonal projection to define the exiting point, i.e., we take the first point of the walk that is outside the domain and then orthogonally project the point onto the boundary of the domain. We consider this projection point as the endpoint for the walk. We define the distribution of this point as discrete exit distribution of SKW with symmetric transition probability, denote it by $\mu_\delta(0, |dz|; D)$. Let $\mu(0, |dz|; D)$ be the harmonic measure for the domain $D$. We expect that $\mu_\delta(0, |dz|; D)$ converges weakly to $\mu(0, |dz|; D)$ as the lattice spacing $\delta$ goes to zero. The difference between these two measures, to first order in the lattice spacing $\delta$, is given by the same form as in the original SKW model,  $\delta c\rho_D(z) |dz|$. Moreover the density function $\rho_D$ is the same as in the usual uniform SKW model and the explicit computation of this function can be found in \cite{tgk}.

From our simulations, the limiting exit distribution of the SKW with the non-uniform but symmetric transition probability as the lattice spacing goes to zero is indeed the harmonic measure.  We also find that the constant $c$ does not depend on the domain, but depends on the transition probability, and the density function $\rho_D$ only depends on the domain $D$. The conjecture in \cite{tgk} also holds in this model. We restate the conjecture here. 

\begin{conjecture}
For each simply connected domain $D$ containing the origin in the plane, there is a function $\rho_D(z)$ on $\partial D$ which only depends on $D$, and for each SKW with symmetric transition probability on the square lattice with lattice spacing $\delta$, there is a constant $c$ which only depends on the transition probability,  such that 
\[\lim_{\delta\to 0} \frac{\mu_\delta(0, |dz|; D)  - \mu(0, |dz|; D)}{\delta} = c\rho_D(z) |dz|\]
Here $|dz|$ is the Lebesgue measure on the boundary with respect to arc length. 
\end{conjecture}

\section{Simulations}
We have generated the data using the same code used in \cite{tgk}. The way we tested the conjecture is to compute the empirical SKW exit distribution and compare it with distribution of the harmonic measure in two particular domains on the unit lattice.  One domain $D_1$ is the circle centered at $(0.3, -0.25)$ with radius 1. Another domain $D_2$ is the horizontal strip of width 1 with the distance from the origin to the top boundary equal to 0.6. Recall that the harmonic measure of the unit disk starting at the origin is length measure on the unit circle. Since the harmonic measure is conformally invariant, if we know the conformal map from the domain to the unit disk, the measure can be easily computed. Finding conformal maps for domains $D_1$ and $D_2$ is an easy exercise. We leave them to the reader.

We are interested in the limit where the lattice spacing  goes to zero. For each domain and case with different transition probability, we ran simulations with $\delta=0.02, 0.01, 0.005$. 

In our simulations, we use cumulative distribution function to display the results. Let $\theta$ represent the polar angle of the orthogonally projected endpoint for the SKW and $H(\theta)$ be the cumulative distribution function (CDF) of the harmonic measure. To test the conjecture, we compared the difference of the empirical cumulative distribution function (ECDF) of the data, denoted as $F(\theta)$, with $H(\theta)$. $10^8$ samples are simulated. 

\begin{figure}[h!]
  \centering
      \includegraphics[width=1\textwidth]{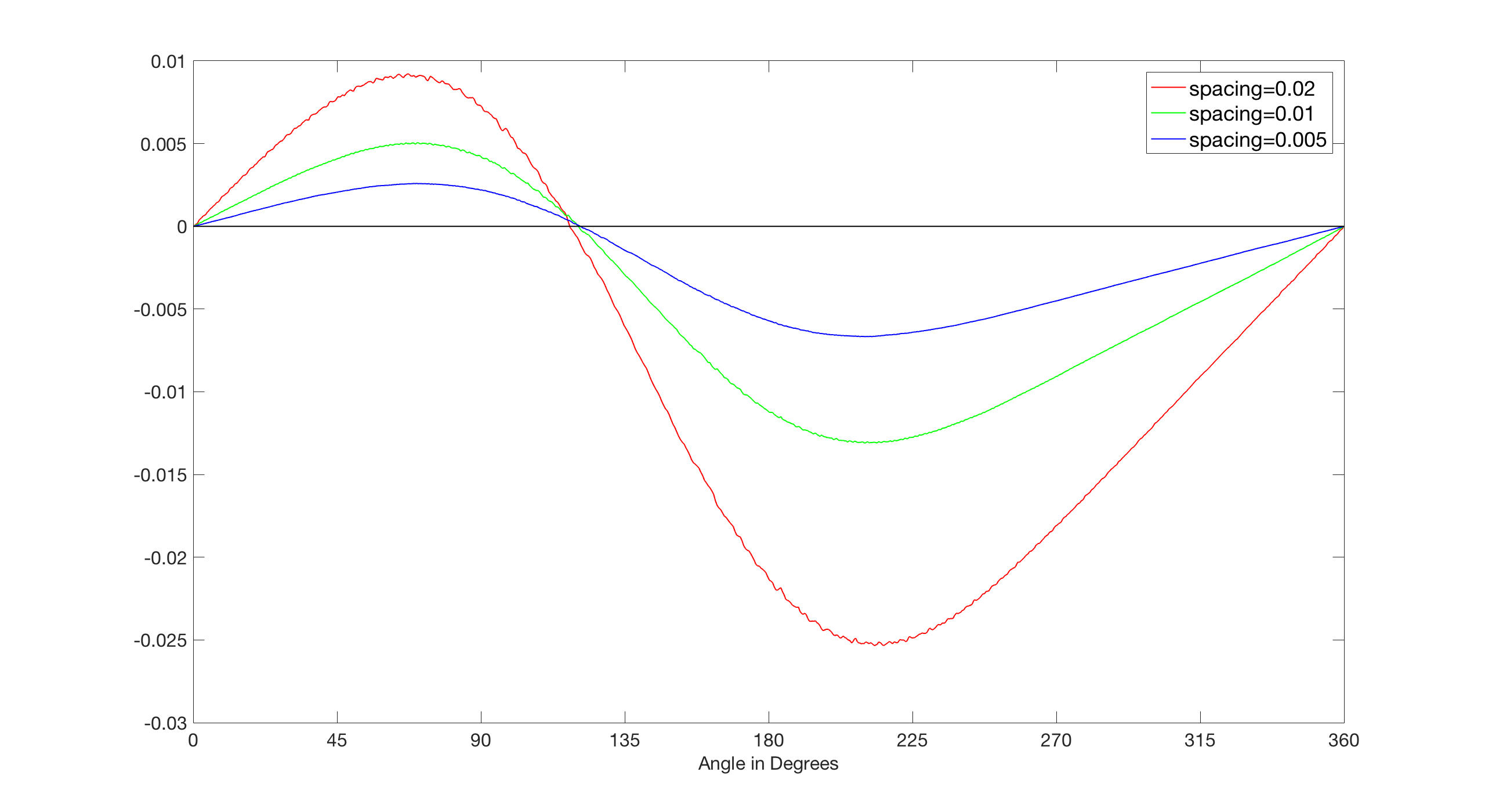}
  \caption{The Difference between $F(\theta)$ and $H(\theta)$ in the case (A) with $a_1$ = 0.9 and $a_2=a_3=0.05$ for domain $D_1$.}
  \label{nc9}
\end{figure}

The first simulation we implemented is a test on the exit distribution of SKW where the front probability $a_1$ in the nblock case (A) has changed and $a_2$ is equal to $a_3$ , while the transition probabilities in (B), (C) and (D) in the figure \ref{transitionprob}, are the same as in the usual uniform SKW model. We have done simulations with $a_1$ = 0.1, 0.3, 0.75, 0.9. Figure \ref{nc9} plots the difference between $F(\theta)$ and $H(\theta)$ for the domain $D_1$ when $a_1$ = 0.9 and $a_2=a_3=0.05$. It shows that the difference between the two functions is going to zero as the lattice spacing $\delta$ goes to zero. This indicates that as the lattice spacing goes to zero, the exit distribution converges to harmonic measure. 

\begin{figure}[h!]
  \centering
      \includegraphics[width=1\textwidth]{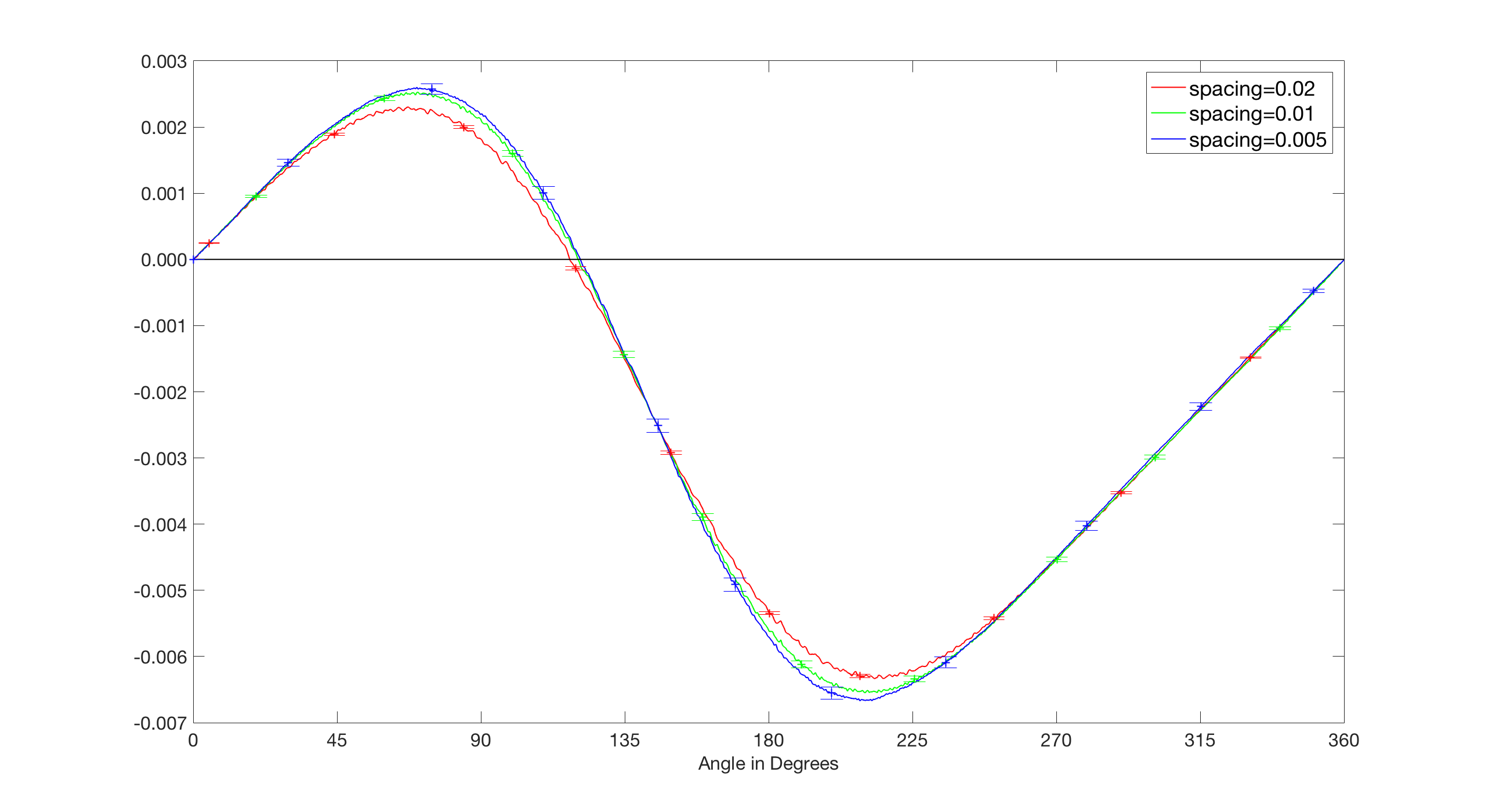}
  \caption{The Rescaled difference  between $F(\theta)$ and $H(\theta)$ in the case (A) with $a_1$ = 0.9 and $a_2=a_3=0.05$ for domain $D_1$.}
  \label{nc9sca}
\end{figure}

\begin{figure}[t]
  \centering
      \includegraphics[width=1\textwidth]{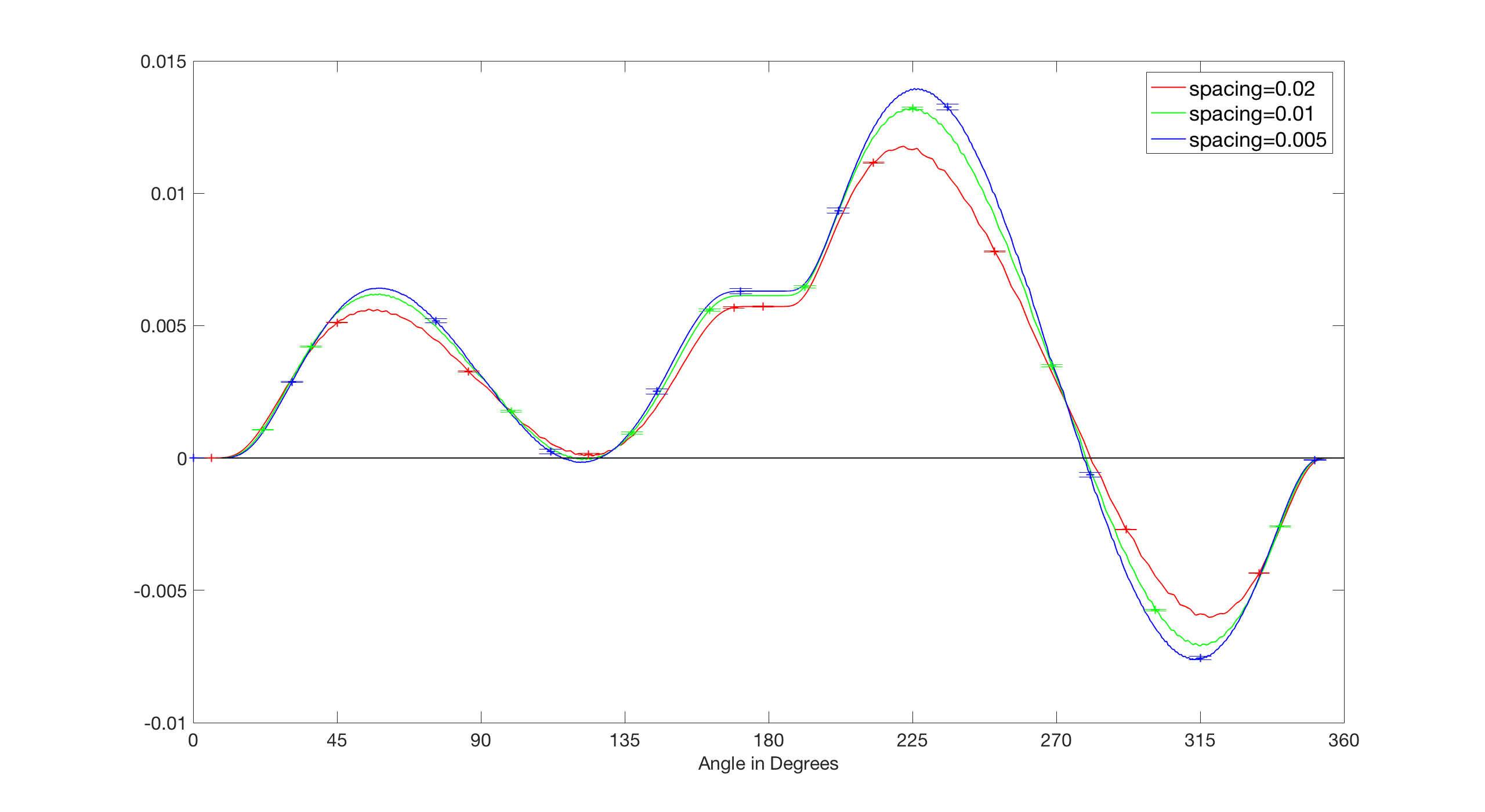}
  \caption{Rescaled difference in the  case (A) with $a_1$ = 0.9 and $a_2=a_3=0.05$ for domain $D_2$.}
  \label{ns9sca}
\end{figure}

In the conjecture, it says that the magnitude of the difference is proportional to the lattice spacing $\delta$. To test this, we fixed the plot for $\delta=0.005$ and vertically compressed the other two plots with $\delta=0.01$ by a factor of 1/2 and $\delta=0.02$ by a factor of 1/4. Figure \ref{nc9sca} shows the rescaled differences between these two functions. As we can see from the figure, three curves collapse nicely into one curve, which supports the conjecture. We also see similar results on the domain $D_2$ for the same case. Figure \ref{ns9sca} shows the result for the rescaled difference with $a_1=0.9$ and $a_2=a_3=0.05$ for the horizontal strip $D_2$. The error bars in both graphs are plus or minus two standard deviations for the statistical errors. There are some places that the curves are not collapsing into one curve even within statistical error. This may come from the fact that the lattice spacing $\delta$ is nonzero, and therefore there are effects beyond first order which may still be significant. 

\begin{figure}[h!]
  \centering
      \includegraphics[width=1\textwidth]{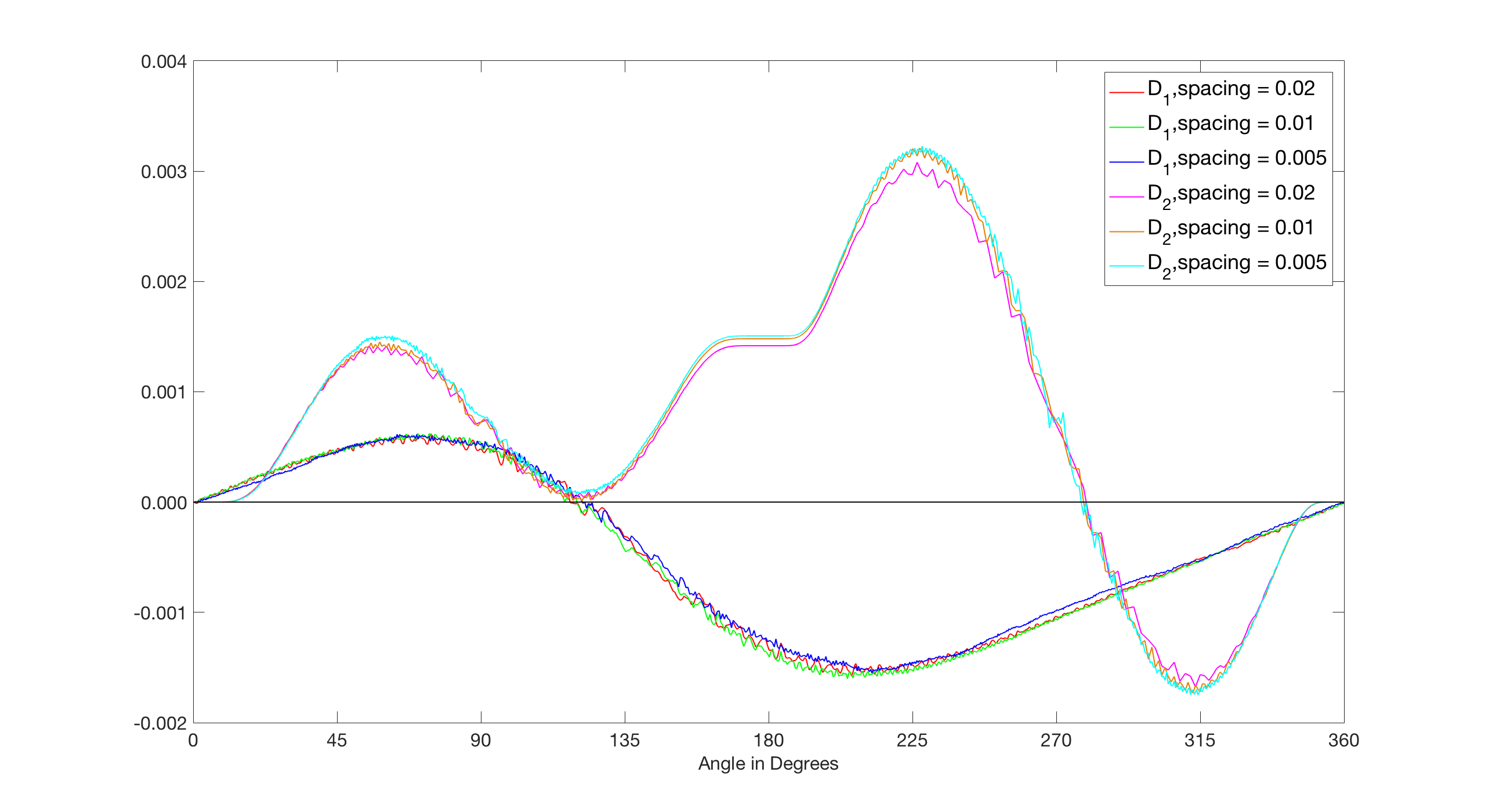}
  \caption{Rescaled difference in the case (B) and (C) with $b_1=c_1=0.1$ for domain $D_1$ and  $D_2$.}
  \label{rl1sca}
\end{figure}

The second part of the test is to run simulations for the scenario where the front probability $b_1$ in the left-blocked case (B) is the same as $c_1$ in the right-blocked case (C),  the transition probabilities in (A) and (D) are kept the same as in the usual uniform SKW model. We ran simulations with $b_1=c_1=0.1, 0.25, 0.33, 0.67, 0.75, 0.9$ in the domain $D_1$. Figure \ref{rl1sca} shows the rescaled differences between the empirical CDF and the CDF of harmonic measure with the $b_1=c_1=0.1$ for the domain $D_1$. Note that the difference is rescaled as we did previously ($\delta=0.005$ is unchanged, $\delta=0.01$ is multiplied by 1/2 and $\delta=0.02$ is multiplied by 1/4). We observe that all three curves collapse into one curve. This again agrees with the conjecture. Figure \ref{rl1sca} also contains three curves for the domain $D_2$; the results are very similar.

\begin{figure}[h!]
  \centering
      \includegraphics[width=1\textwidth]{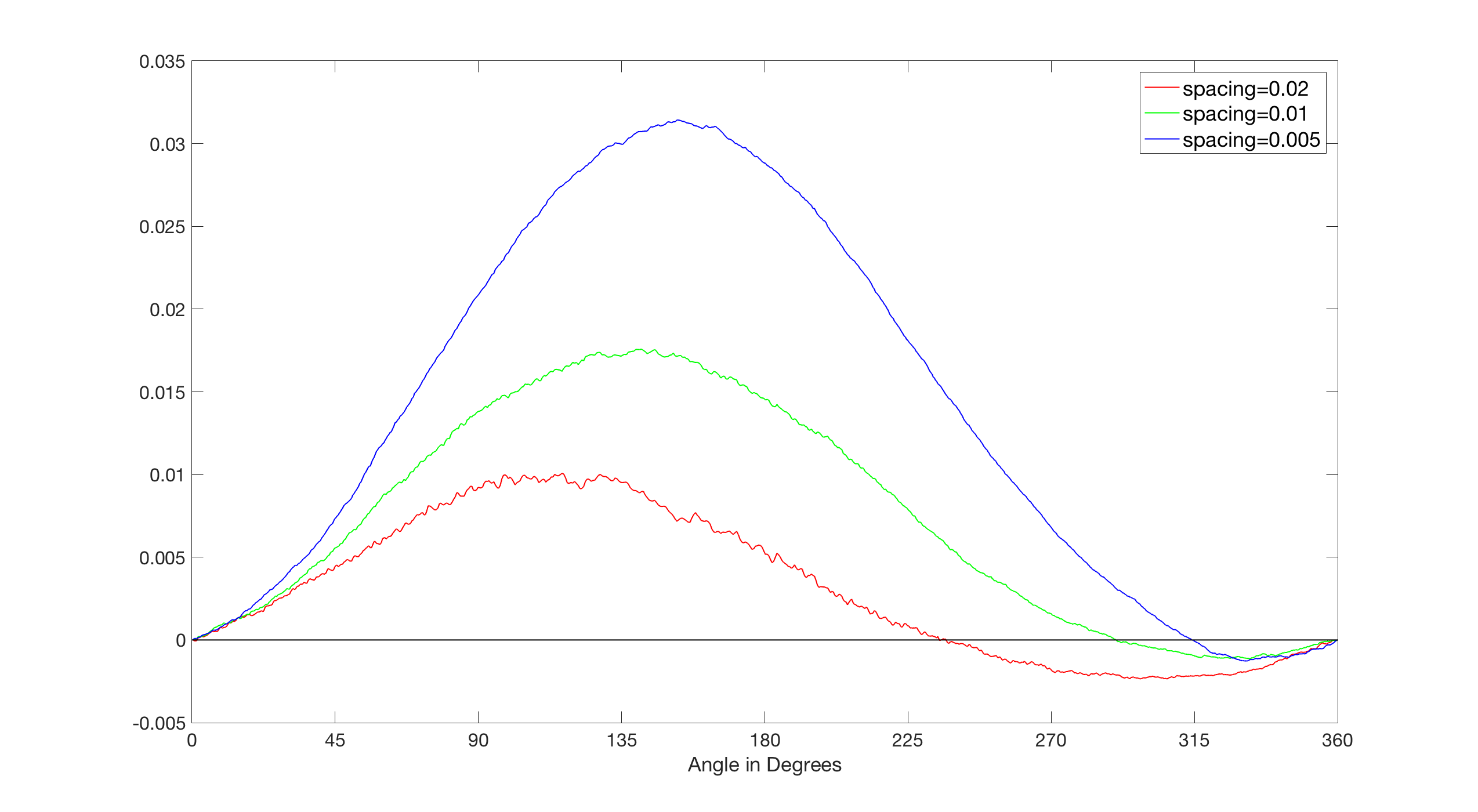}
  \caption{Difference with front probability $b_1$ = 0.55 in the case (B) for domain $D_1$.}
  \label{lc49}
\end{figure}

\begin{figure}[h!]
  \centering
      \includegraphics[width=1\textwidth]{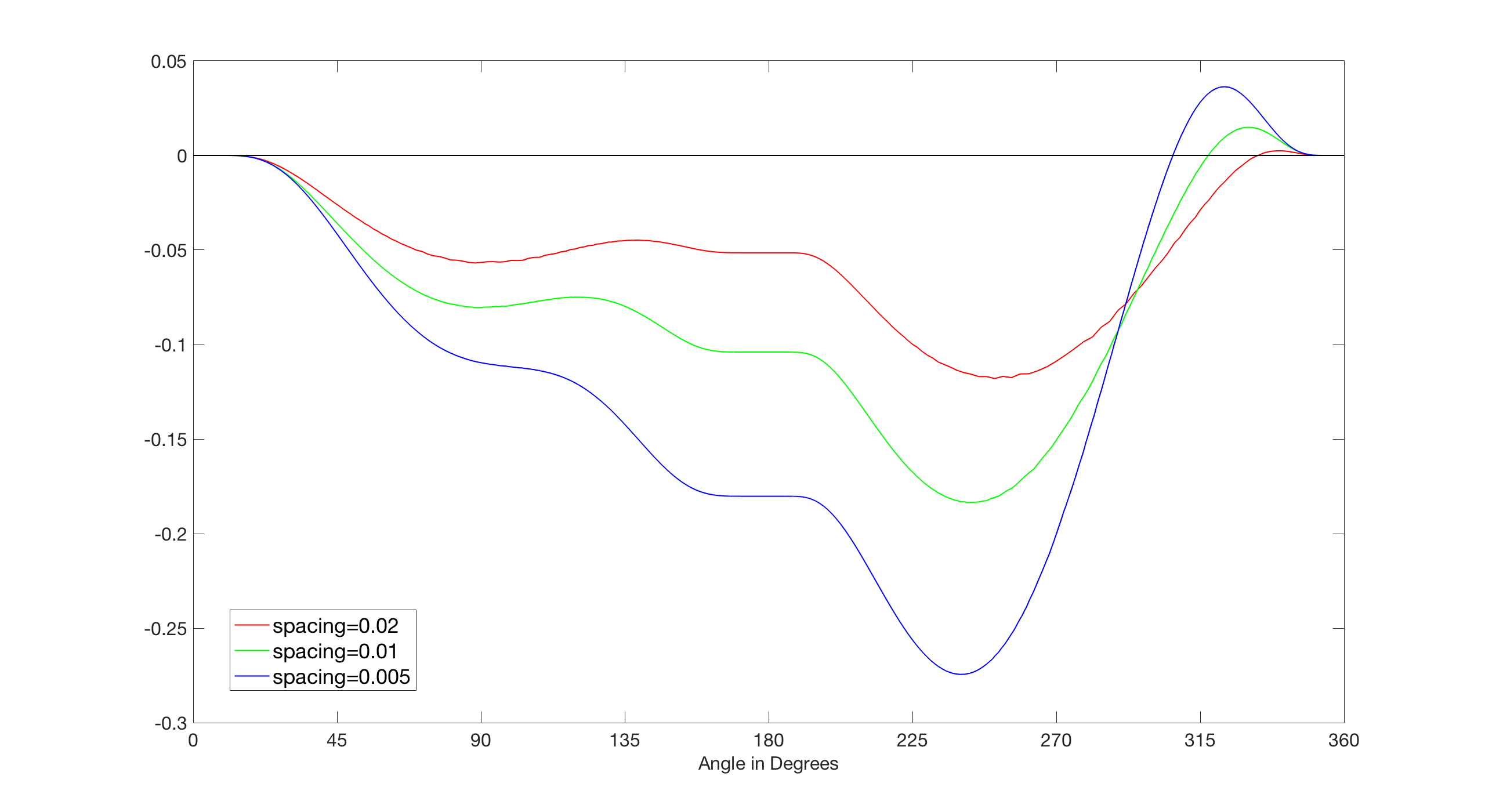}
  \caption{Difference with $a_1=a_2 = 0.3$ and $a_3=0.4$ in the  case (A) for domain $D_2$.}
  \label{ns4}
\end{figure}

From the simulations, we see that the exit distribution of symmetric SKW converges to harmonic measure as the lattice spacing goes to zero. It is interesting to know whether this result holds when the transition probability is not symmetric. We have run several asymmetric symmetric transition probabilities simulations on the nblock case, the left-blocked case and the right-blocked case. Even with a small probability change compared to the transition probability in the usual uniform SKW model,  for example, changing the front probability $b_1$ to 0.55 in the left-blocked case (B) and keeping the rest of the transition probability the same as in the original model for the domain $D_1$, the exit distribution of SKW does not converge to harmonic measure. Figure \ref{lc49} shows the difference between $F(\theta)$ and $H(\theta)$ with $b_1$ = 0.55 in the left-blocked case (B) for the domain $D_1$. In figure \ref{ns4}, we plot the difference function with $a_1=a_2= 0.3$ and $a_3=0.4$ in the nblocked case (A) for the domain $D_2$. We see that in these two plots the difference is not going to zero, and so the exit distribution of the asymmetric SKW  does not converge to harmonic measure as the lattice spacing goes to zero. 

\begin{table}[t]
\begin{center}
\begin{tabular}{c|c|c|}
\cline{2-3}
& $D_1$  & $D_2$\\ \hline
 \multicolumn{1}{ |c| }{ $a_1=0.1$}&    0.8043$\pm$0.0411  & 0.7732 $\pm$0.1061\\ \cline{1-3}
\multicolumn{1}{ |c| }{ $a_1=0.9$}& 6.2874$\pm$0.1729 & 6.3463$\pm$0.1872 \\ \hline
 \multicolumn{1}{ |c| }{ $b_1 = c_1 =0.1$}&1.4055$\pm$0.1797&1.4948 $\pm$0.0658\\ \cline{1-3}
\multicolumn{1}{ |c| }{ $b_1 = c_1 =0.9$}&2.2932$\pm$0.1494 & 2.1596$\pm$0.1523\\ \cline{1-3}
\end{tabular}
\caption{The ratio of the $L^1$ norm of the simulation difference function to the $L^1$ norm of $\rho_D$ with front probability in each case  indicated.}
  \label{table1}
  \end{center}
\end{table}

The third part of this section is to study the density function $\rho_D$. In the conjecture, it says that $\rho_D$ does not depend on the transition probability, and only depends on the domain. As we can see from figure \ref{nc9} to \ref{rl1sca}, there are only two shapes of curves. Comparing the curves for the domain $D_1$ in figure \ref{nc9sca} and \ref{rl1sca}, they are quite similar and only with different vertical scales. Likewise, the curves in figure \ref{ns9sca} and \ref{rl1sca} for the domain $D_2$ are also similar. This indicates the the function $\rho_D$ has a universal shape for the same domain except for an overall different vertical scale. 

Finally we test the dependence of the constant $c$ in the conjecture. The conjecture says that the constant $c$ only depends on the transition probability, not the domain. In order to test this, we compute the ratio of the $L^1$ norm of the simulation difference function to the $L^1$ norm of $\rho_D$. For the comparison, we consider the simulations with the lattice spacing $\delta=0.005$. In the simulations, the domain $D_1$ in table \ref{table1} refers to the circle centered at (60, $-$50) with radius 200, and $D_2$ is the horizontal strip of width 200 with distance from the origin to the top boundary equal to 120. Table \ref{table1}  shows the ratio for two domains with 95\% confidence intervals. If the conjecture is correct, then the values in the table should be the same for the different domains and same transition probability. We found that the values of each row in table \ref{table1}  are quite close within statistical error, which supports the conjecture.

\section{Conclusions}
In this paper, we have studied the SKW with non-uniform but symmetric or asymmetric transition probability with averaging over rotations of the lattice. We see from the simulations that the exit distribution of walk with symmetric  transition probability does converge weakly to harmonic measure as the lattice spacing goes to zero. The simulations also support the conjecture that the difference of the exit distribution of this model and the harmonic measure, to first order in the lattice spacing $\delta$, is given by $\delta c\rho_D(z) |dz|$. Moreover, the constant $c$ only depends on the transition probability, and the function $\rho_D$ only depends on the domain. For the SKW with asymmetric transition probability, the exit distribution does not converge. 

It would be interesting to do the same test on the triangular lattice for the SKW model with symmetric transition probability. We would expect that the limiting distribution of SKW with symmetric transition probability as the lattice spacing goes to zero is harmonic measure. The first order correction on the triangular lattice, up to an overall constant, is the same as on the square lattice. The constant would depend on the transition probability and the lattice, but not on the domain. On the hexagonal lattice, there are only two possible nearest neighbors for each site. If the SKW model is symmetric on the hexagonal lattice, then the probabilities of these two nearest neighbors have to be the same. 

\vspace{0.5cm}
\text{\bf{Acknowledgments}} This research was supported in part by NSF grant DMS-1500850. The numerical computations reported here were performed at the UA Research Computing High Performance Computing (HPC) and High Throughput Computing (HTC) at the University of Arizona. The author would like to thank Tom Kennedy for providing the author the  guidance and support for the research. 
\bibliographystyle{amsplain}

\end{document}